\documentclass[11pt,reqno]{amsart}
\usepackage{amsmath,amssymb,amsthm,upref,graphicx,mathrsfs}

\usepackage{color,xcolor}
\usepackage[
  colorlinks=true,
  linkcolor=blue,
  citecolor=blue,
  urlcolor=blue]{hyperref}

\newtheorem{thm}{Theorem}[section]
\newtheorem{cor}[thm]{Corollary}

\newtheorem{prop}[thm]{Proposition}

\newtheorem{exaple}[thm]{Example}
\numberwithin{equation}{section}

\begin{document}

\leftline{ \scriptsize}

\vspace{1.3 cm}
\title
{Generalized matrix functions on a linear sum of permutation matrices and their cousins}
\author{Ratisiri Sanguanwong $^{a}$ and Kijti Rodtes$^{b,\ast}$ }
\thanks{{\scriptsize
$^{\ast}$ Corresponding author.\\
E-mail addresses: r.sanguanwong@gmail.com (Ratsiri Sanguanwong), kijtir@nu.ac.th (Kijti Rodtes).\\
$^{A}$ Department of Mathematics, Faculty of Science, Naresuan University, Phitsanulok 65000, Thailand.\\
$^{B}$ Department of Mathematics, Faculty of Science, Naresuan University, and Research Center for Academic Excellent in Mathematics, Phitsanulok 65000, Thailand.\\}}
\hskip -0.4 true cm

\maketitle


\begin{abstract}A generalized matrix function is a generalization of determinant and permanent function. In this paper, we introduced the formula for the value of a generalized matrix function of a linear sum of permutation matrices. We show that a linear sum of permutation matrices satisfies the permanent dominance conjecture. Finally, we apply the result to some cousins of permutation matrices.
\end{abstract}

\vskip 0.2 true cm


\pagestyle{myheadings}
\markboth{\rightline {\scriptsize Ratsiri Sanguanwong and  Kijti Rodtes}}
         {\leftline{\scriptsize }}
\bigskip
\bigskip


\vskip 0.4 true cm

\section{Introduction}

A generalized matrix function is a function $d_\chi^G:M_n(\mathbb{C}) \rightarrow \mathbb{C}$, constructed by a subgroup $G$ of $S_n$ and a character $\chi$ of $G$. This function is a generalization of determinant and permanent function. Most papers relating to generalized matrix functions regarded to positive semidefinite (PSD) matrices. In 1965, M. Marcus and H. Minc established a relation between a generalized matrix function and singular values of a given matrix \cite{MM1965}. In 1976, R. Merris obtained some relations among generalized matrix functions on PSD matrices \cite{RM1976}. In 2017, S. Huang et. al. investigated inequalities of generalized matrix functions on PSD matrices \cite{SH2017}.

For any two matrices $A, B$, a relation between $\det(A + B)$ and $\det(A)+\det(B)$ is obtained by Cauchy-Binet formula; namely, by applying the formula to matrices $\begin{bmatrix}
A & I
\end{bmatrix}$ and $\begin{bmatrix}
I\\B
\end{bmatrix}$, we derive that (cf. \cite{RMerrisB})
\begin{equation}\label{eq00}
\det(A + B) = \sum\limits_{k = 0}^n\sum\limits_{\alpha,\beta\in\Omega_{k,n}}(-1)^{r(\alpha)+r(\beta)}\det(A[\alpha\mid\beta])\det(B(\alpha\mid\beta)),
\end{equation}
where $r(\alpha) = \sum\limits_{i=1}^k \alpha(i)$, matrix $A[\alpha\mid\beta] (\hbox{resp. }A(\alpha\mid\beta))$ is a submatrix of $A$ obtained by keeping (resp. deleting) row and column corresponding to $\alpha$ and $\beta$, respectively. Here, $\Omega_{m,n} = \{(\omega_1, \dots, \omega_m)\mid 1\leq\omega_1<\dots<\omega_n\leq n\}$. However, to calculate the value of $\det(A + B)$ by this formula, we need to construct $\sum\limits_{k = 0}^n \binom{n}{k}^2$ submatrices, even for the case $aP_\theta + bP_\tau$, a linear sum of permutation matrices corresponding to $\sigma$ and $\tau$.
For the case of PSD matrices $A, B$, we have that \cite{RMerrisB}
\begin{equation}\label{eq0}
d_\chi^G(A + B) \geq d_\chi^G(A) + d_\chi^G(B),
\end{equation}
where $\chi$ is an irreducible character. This provide only a lower bound for $d_\chi^G(A + B)$, but not the exact value.

Permutation matrices have many applications and properties \cite{RMerrisB, FZhangB}. A concept of these matrices can be extended to generalized permutation matrices. In \cite{RS2019}, the authors worked with symmetric matrices constructed by permutations. These matrices can be viewed as a sum of two generalized permutation matrices, denoted by $S_\sigma$. There exists $\sigma\in S_n$ such that $S_\sigma = P_\sigma + P_{\sigma^{-1}}$ while $d_\chi^G(S_\sigma) < d_\chi^G(P_\sigma)+d_\chi^G(P_\sigma^{-1})$. This means that the equality (\ref{eq0}) does not hold when $A,B$ are generalized permutation matrices. We refer to such symmetric matrices and generalized permutation matrices as cousins of permutation matrices.

In this paper, we pay attention to permutation matrices and their cousins. The value of generalized matrix function of each permutation matrix is well-known. Our work is the investigation of the formula for a generalized matrix function of a linear sum of two permutation matrices. A simple formula is obtained for some special cases. A necessary and sufficient condition for definiteness of these matrices is established. A result relating to a singular value of these matrices is remarked. We also show that a linear sum of permutation matrices satisfies the permanent dominance conjecture. Moreover, we provide the formula to some cousins of the permutation matrices.

\section{Preliminary}
Let $V$ be an $n$-dimensional complex vector space. Then  $GL(V)$, the set of all bijective linear transformations of $V$, is a group under composition of functions. For each finite group $G$, there exists a homomorphism $\rho:G \rightarrow GL(V)$. A \textit{character} $\chi$ of $G$ is a mapping from $G$ into $\mathbb{C}$ defined by $\chi(g) = tr(\rho(g))$ for all $g \in G$. A character $\chi$ of $G$ is said to be \textit{linear} if it is a homomorphism.

For each $n \in \mathbb{N}$, denote $[n] = \{1,2,\dots,n\}$. A \textit{symmetric group} of degree $n$, denoted by $S_n$, is a set of all bijections on $[n]$. An element of $S_n$ is called a \textit{permutation}. A permutation $\sigma\in S_n$ is called a \textit{cycle} if there exists a subset $S$ of $[n]$ such that $\operatorname{Fix}(\sigma) = [n]\setminus S$ and for each $i,j \in S$, there exists $k \in \mathbb{N}$ satisfying $j = \sigma^k(i)$.  A cycle is called $k$-cycle if it is an element of order $k$. For each $\sigma \in S_n$, it can be written uniquely as a product of disjoint cycles (not including a $1$-cycle) $C_1\cdots C_r$. Here, we say that $\sigma$ has a \textit{cycle structure} $[l_1, \dots, l_r, 1^F]$, where $C_i$ is an $l_i$-cycle and $F$ is a number of fixed points of $\sigma$. In this paper, the product of disjoint cycles not including $1$-cycle is called the \textit{disjoint cycles expression}. For each subgroup $G$ of $S_n$ and a character $\chi$ of $G$, a \textit{generalized matrix function} is a function $d_\chi^G:M_n(\mathbb{C})\rightarrow\mathbb{C}$ defined by
\begin{equation*}
\hbox{$d_\chi^G(A) = \sum\limits_{\sigma\in G}\chi(\sigma)\prod\limits_{i=1}^n A_{i \, \sigma(i)}$,}
\end{equation*}
where $A_{i \, j}$ denote an $(i, j)$-entry of matrix $A$.
By considering this formula, determinant and permanent function are special cases of generalized matrix functions. For $\theta \in S_n$, a \textit{permutation matrix} corresponding to $\theta$, denote by $P_\theta$, is a square matrix such that $(P_\theta)_{\theta(j) \, j}=1$ for each $j=1,\dots,n$ and the other entires are zero. In other words, we say that
\begin{equation}\label{pm}
(P_\theta)_{i \, j} = \left\{
\begin{array}{ll}
1, & \hbox{if $\theta^{-1}(i)=j$,} \\
0, & \hbox{otherwise.}
\end{array}
\right.
\end{equation}
In this paper, we consider a permutation matrix as (\ref{pm}). Let $D$ be an invertible diagonal matrix and $\sigma \in S_n$. A matrix $G = DP_\sigma$ is called a \textit{generalized permutation matrix}.

For any matrix $A$, denote a \textit{conjugate transpose} of $A$ by $A^\ast$. A complex square matrix $A$ is called a \textit{hermitian matrix} if $A = A^{\ast}$. A complex number $\lambda$ is said to be an \textit{eigenvalue} of an $n \times n$ matrix $A$ if $\det(\lambda I - A) = 0$. An eigenvalue of hermitian matrix is a real number. A hermitian matrix $A$ is said to be \textit{positive semidefinite}, briefly \textit{PSD}, if all eigenvalues of $A$ are non-negative. The permanent dominance conjecture states that, for each PSD matrix $A$ and a generalized matrix function $d_\chi^G$, \cite{RMerrisB},
\begin{equation*}
\hbox{$\frac{\displaystyle 1}{\displaystyle\chi(id)}d_\chi^G(A)\leq \operatorname{per}(A)$.}
\end{equation*}

A non-negative square root of an eigenvalue of $A^\ast A$ is called a \textit{singular value} of $A$. For any square matrix $A$ with singular values $\alpha_1 , \dots ,\alpha_n$, we have, \cite{MM1965},
\begin{equation}\label{sveq}
\hbox{$\lvert d_\chi^G(A) \rvert\leq \frac{\displaystyle 1}{\displaystyle n}\sum\limits_{i=1}^n\alpha_i^{2n},$}
\end{equation}
where $\chi$ is a linear character of a subgroup $G$ of $S_n$.

Let $V$ be an inner product space. Let $V^{(n)}$ denote the $n$-folds tensor product $V \otimes \cdots \otimes V$. For $\sigma \in S_n$, let $P(\sigma): V^{(n)}\rightarrow V^{(n)}$ be the permutation operator with
\begin{equation*}
P(\sigma)(x_1 \otimes \cdots \otimes x_n) = x_{\sigma^{-1}(1)}\otimes \cdots \otimes x_{\sigma^{-1}(n)}.
\end{equation*}
Define a symmetry operator $T:V^{(n)} \rightarrow V^{(n)}$ by
\begin{equation*}
\hbox{$T = \sum\limits_{\sigma \in G}\chi(\sigma)P(\sigma)$}.
\end{equation*}
We denote $x_1 \ast \cdots \ast x_n := T(x_1 \otimes \cdots \otimes x_n)$. Let $x_1, \dots , x_n$ and $y_1, \dots, y_n$ be a vector in $V$. If $A = (\langle x_i, y_j \rangle)$, then, \cite{MM1965},
\begin{equation}\label{eq21}
\langle x_1 \ast \cdots x_n, y_1 \ast \cdots \ast y_n \rangle = \lvert G \rvert d_\chi^G(A).
\end{equation}
\section{Generalized matrix functions on a linear sum of two permutation matrices}

\quad Let $\theta,\tau \in S_n$. We will evaluate the value of $d_\chi^G(aP_\theta + bP_\tau)$, where $a,b \in \mathbb{C}$. By the definition of $d_\chi^G$, we have
\begin{eqnarray*}
d_\chi^G(aP_\theta + bP_\tau)
& = &\sum\limits_{\sigma \in G}\chi(\sigma)\prod\limits_{i=1}^n(a(P_\theta)_{i \, \sigma(i)}+b(P_\tau)_{i \, \sigma(i)})\label{eq01}
\end{eqnarray*}

Moreover,
\begin{equation}\label{eq05}
\hbox{$\prod\limits_{i=1}^n(a(P_\theta)_{i \, \sigma(i)}+b(P_\tau)_{i \, \sigma(i)}) = 0 $ if $\exists i \in [n], \sigma(i) \neq \theta^{-1}(i)$ and $\sigma(i) \neq \tau^{-1}(i)$.}
\end{equation}

Denote $X(\theta, \tau) := \{\sigma\in S_n \mid \hbox{for each } i \in [n], \sigma(i) = \theta(i) \hbox{ or } \sigma(i) = \tau(i)\}$. By this definition, we have $X(\theta^{-1},\tau^{-1}) = \{\sigma^{-1} \mid \sigma \in X(\theta,\tau)\}$. Then, by (\ref{eq05}),
\begin{equation*}
d_\chi^G(aP_\theta + bP_\tau) = \sum\limits_{\sigma \in X(\theta, \tau)\cap G}\chi(\sigma^{-1})\prod\limits_{i=1}^n(a(P_\theta)_{i \, \sigma^{-1}(i)}+b(P_\tau)_{i \, \sigma^{-1}(i)}).
\end{equation*}

For each $\sigma \in X(\theta, \tau)$, $\sigma(i) = \theta(i) = \tau(i)$ if and only if $i \in \operatorname{Fix}(\theta^{-1}\tau)$. For $\pi \in S_n$ and $\sigma \in X(\theta,\tau)$, we denote $I^{\sigma}_\pi=\{i \in [n] \mid \sigma(i)=\pi(i)\}$. Then $[n]$ is partitioned into the disjoint union as follow:
\begin{equation}\label{partition}
(I^{\sigma}_\theta\setminus I^{\sigma}_\tau) \cup (I^{\sigma}_\tau\setminus I^{\sigma}_\theta) \cup (I^{\sigma}_\theta\cap I^{\sigma}_\tau) = [n].
\end{equation}
Here, we have $I^{\sigma}_\theta \cap I^{\sigma}_\tau = \operatorname{Fix}(\theta^{-1}\tau)$, for each $\sigma \in X(\theta,\tau)$.

\begin{thm}\label{thm01}
	Let $\theta,\tau\in S_n$. Suppose that $\theta^{-1}\tau = C_1 C_2 \cdots C_r$ is the disjoint cycles expression. Then $X(\theta, \tau)=\{\theta\prod\limits_{i \in I}C_i \mid I\subseteq [r]\}$.
\end{thm}
\begin{proof}
	Let $\sigma \in X(\theta, \tau)$. We claim that, for each $j \in [r]$, if $C_j$ contains an element of $I^{\sigma}_\theta$, then $C_j$ does not contain any element of $[n]\setminus I^{\sigma}_\theta$. Suppose that $C_j$ contain $k, l$, where $k \in [n]\setminus I^{\sigma}_\theta$ and $l \in I^{\sigma}_\theta$. Then $\sigma(k)=\tau(k)$ and $\sigma(l)=\theta(l)$. Without loss of generality, we may assume that $C_j(k)=l$. This implies that
	\begin{equation*}
	\theta^{-1}\sigma(k) = \theta^{-1}\tau(k)=C_j(k)=l.
	\end{equation*}
	Thus $\sigma(k) = \theta(l)=\sigma(l)$. This is a contradiction. Now, we have the claim. Suppose that $C_1,\ldots, C_t$ are all cycles containing elements of $I^{\sigma}_\theta$. Let $I = [r]\setminus[t]$. We can evaluate that
	\begin{equation*}
	\hbox{$\sigma(x) = \theta(x) = \theta\prod\limits_{i \in I}C_i(x)$}
	\end{equation*}
	if $x$ is not contained in $C_i$ for each $i \in I$. On the other hand,
	\begin{equation*}
	\hbox{$\sigma(x) = \tau(x) = \theta\prod\limits_{i \in I}C_i(x)$}
	\end{equation*}
	if $x$ is contained in $C_i$ for some $i \in I$.
 	This means that that $\sigma = \theta\prod\limits_{i \in I}C_i$. Thus $X(\theta, \tau) \subseteq \{\theta\prod\limits_{i \in I}C_i \mid I\subseteq [r]\}$. Conversely, let $\sigma = \theta\prod\limits_{i \in I}C_i$ for some $I\subseteq [r]$. By a similar reasoning as above, we have $\sigma(x) = \theta(x)$ if $x$ is not contained in $C_i$ for every $i \in I$, otherwise, $\sigma(x) = \tau(x)$. Thus $\sigma \in X(\theta,\tau)$.
\end{proof}

Next, we shall evaluate the value of $\prod\limits_{i=1}^n(a(P_\theta)_{i \, \sigma^{-1}(i)}+b(P_\tau)_{i \, \sigma^{-1}(i)})$, where $\sigma \in X(\theta, \tau)$. For each $I \subseteq [n]$, we can see that
\begin{equation}
\hbox{$\prod\limits_{i\in I}a(P_\theta)_{i \, \sigma(i)}\prod\limits_{i \not\in I}b(P_\tau)_{i \, \sigma(i)}= 0 $ if $\exists i \in I, \sigma(i) \neq \theta^{-1}(i)$ or $\exists i \not\in I, \sigma(i) \neq \tau^{-1}(i)$}.\label{eq04}
\end{equation} For each subset $A$ of $\operatorname{Fix}(\theta^{-1}\tau)$, we denote $I^{\sigma}_A = (I^{\sigma}_\theta\setminus I^{\sigma}_\tau)\cup A$. If $I$ is a subset of $[n]$ such that $\prod\limits_{i\in I}a(P_\theta)_{i \, \sigma(i)}\prod\limits_{i \not\in I}b(P_\tau)_{i \, \sigma(i)}\neq 0$, then, by (\ref{eq04}),
\begin{equation*}
\sigma(i) = \left\{
\begin{array}{ll}
\theta^{-1}(i), & \hbox{if $i \in I$,} \\
\tau^{-1}(i), & \hbox{otherwise.}
\end{array}
\right.
\end{equation*}
So, $(I^{\sigma}_\theta\setminus I^{\sigma}_\tau) \subseteq I$ and $(I^{\sigma}_\tau\setminus I^{\sigma}_\theta) \subseteq [n]\setminus I$. By (\ref{partition}),
\begin{equation*}
(I^{\sigma}_\theta\setminus I^{\sigma}_\tau) \subseteq I \subseteq (I^{\sigma}_\theta\setminus I^{\sigma}_\tau) \cup \operatorname{Fix}(\theta^{-1}\tau).
\end{equation*}
Thus, $I = I^{\sigma}_A$ for some $A \subseteq \operatorname{Fix}(\theta^{-1}\tau)$ if $\prod\limits_{i\in I}a(P_\theta)_{i \, \sigma(i)}\prod\limits_{i \not\in I}b(P_\tau)_{i \, \sigma(i)}\neq 0$. Suppose that $\theta^{-1}\tau = C_1 C_2 \cdots C_r$ is the disjoint cycles expression. By Theorem \ref{thm01}, there exists $J\subseteq[r]$ such that $\sigma = \theta\prod\limits_{i\ \in J}C_i$. Put
\begin{equation}\label{tF}
\hbox{$t_\sigma = n-\lvert \operatorname{Fix}(\theta^{-1}\sigma) \rvert$ and $F = \lvert \operatorname{Fix}(\theta^{-1}\tau) \rvert$.}
\end{equation}
It is obvious that $t_{\sigma^{-1}}=t_\sigma$. Then, for each $\sigma \in X(\theta,\tau)$,
\begin{eqnarray*}
\prod\limits_{i=1}^n(aP_\theta + bP_\tau)_{i \, \sigma(i)} & = & \sum\limits_{I \subseteq [n]}(\prod\limits_{i\in I}a(P_\theta)_{i \, \sigma(i)}\prod\limits_{i \not\in I}b(P_\tau)_{i \, \sigma(i)})\\
 & = & \sum\limits_{A \subseteq \operatorname{Fix}(\theta^{-1}\tau)}(\prod\limits_{i\in I^{\sigma}_A}a(P_\theta)_{i \, \sigma(i)}\prod\limits_{i \not\in I^{\sigma}_A}b(P_\tau)_{i \, \sigma(i)})\\
 & = & \sum\limits_{A \subseteq \operatorname{Fix}(\theta^{-1}\tau)}a^{\lvert I^{\sigma}_A \rvert} b^{n - \lvert I^{\sigma}_A \rvert}\\
 & = & \sum\limits_{A \subseteq \operatorname{Fix}(\theta^{-1}\tau)}a^{n - t_\sigma - \lvert A \rvert} b^{t_\sigma + \lvert A \rvert}\\
 & = & \sum\limits_{k=0}^F \binom{F}{k} a^{n - t_\sigma - k} b^{t_\sigma + k}\\
 & = & a^{n - t_\sigma - F} b^{t_\sigma}(a+b)^F.
\end{eqnarray*}
For a character $\chi$ of $G$, we have $\overline{\chi}(\sigma) = \chi(\sigma^{-1})$ for each $\sigma \in G$. The following consequence is obtained immediately.
\begin{thm}\label{form1}
	Let $\theta,\tau \in S_n$ and $a,b \in \mathbb{C}$. Then
	\begin{equation*}
	d_\chi^G (aP_\theta + bP_\tau) = (a+b)^F\sum\limits_{\sigma\in X(\theta,\tau) \cap G}\overline{\chi}(\sigma)a^{n - t_\sigma - F} b^{t_\sigma},
	\end{equation*}
	where $t_\sigma$ and $F$ are defined as (\ref{tF}).
\end{thm}

By using Theorem \ref{form1} instead of (\ref{eq00}), we only need to find $X(\theta,\tau)$ which contains $2^r$ elements, where $r$ is the number of cycles in the disjoint cycles expression of $\theta^{-1}\tau$.
\begin{exaple}
	Let $\theta = (153)(26)$ and $\tau = (246)$ be permutations in $S_6$. By direct calculation, $\theta^{-1}\tau = (135)(24)$. Then $X(\theta,\tau) = \{\theta,\tau, (26), (153)(246)\}$. Then, by Theorem \ref{form1},
	\begin{equation*}
	\det(\begin{bmatrix}
	-i & 0 & 2 & 0 & 0 & 0\\
	0 & 0 & 0 & 0 & 0 & 2-i\\
	0 & 0 & -i & 0 & 2 & 0\\
	0 & -i & 0 & 2 & 0 & 0\\
	2 & 0 & 0 & 0 & -i & 0\\
	0 & 2 & 0 & -i & 0 & 0
	\end{bmatrix}) = (2-i)(-32 -i -4i - 8) = -85+30i.
	\end{equation*}
	Let $G = \{\sigma\in S_6 \mid 1,3,5 \in \operatorname{Fix}(\sigma)\}$ and $\chi$ a character of $G$. Here, $X(\theta,\tau) \cap G = \{\tau,(26)\}$. For each $a, b \in \mathbb{C}$, we have
	\begin{equation*}
	d_\chi^G(\begin{bmatrix}
	b & 0 & a & 0 & 0 & 0\\
	0 & 0 & 0 & 0 & 0 & a+b\\
	0 & 0 & b & 0 & a & 0\\
	0 & b & 0 & a & 0 & 0\\
	a & 0 & 0 & 0 & b & 0\\
	0 & a & 0 & b & 0 & 0
	\end{bmatrix}) =(a+b)(\overline{\chi}(\tau)b^5 + \overline{\chi}((26))a^2b^3).
	\end{equation*}
\end{exaple}
By applying this formula, for determinant and permanent functions, we only need to consider the cycles structure of $\theta^{-1}\tau$. Some cycle structures yield more simple results.
\begin{cor}
For $\theta,\tau \in S_n$ and $a,b\in \mathbb{C}$, we have
\begin{equation*}
\det(aP_\theta + bP_\tau) = \operatorname{sgn}(\theta) \sum\limits_{I \subseteq [r]}(-1)^{\lvert I \rvert + \sum\limits_{i\in I}l_i} a^{n-\sum\limits_{i\in I}l_i-F} b^{\sum\limits_{i\in I}l_i} (a+b)^F
\end{equation*}
and
\begin{equation*}
\operatorname{per}(aP_\theta + bP_\tau) = \sum\limits_{I \subseteq [r]}a^{n-\sum\limits_{i\in I}l_i-F} b^{\sum\limits_{i\in I}l_i} (a+b)^F,
\end{equation*}
	 where $[l_1, l_2, ..., l_r, 1^F]$ is the cycle structure of $\theta^{-1}\tau$.
\end{cor}

\begin{cor}
	For $\theta,\tau \in S_n$ and $a,b\in \mathbb{C}$, if $\theta^{-1}\tau$ is an $n$-cycle, then
	\begin{equation*}
	\det(aP_\theta + bP_\tau) = \operatorname{sgn}(\theta) (a^n + (-b)^n)=\det(aP_\theta)+\det(bP_\tau)
	\end{equation*}
	and
	\begin{equation*}
	\operatorname{per}(aP_\theta + bP_\tau) = a^n + b^n = \operatorname{per}(aP_\theta)+\operatorname{per}(bP_\tau).
	\end{equation*}
\end{cor}

\begin{cor}
	For $\theta,\tau \in S_n$ and $a,b\in \mathbb{C}$, if all $r$ cycles in the disjoint cycles expression of $\theta^{-1}\tau$ have the same length $l$, then
	\begin{equation*}
	\det(aP_\theta + bP_\tau) = \operatorname{sgn}(\theta) (a^l + (-1)^{l+1}b^l)^r (a+b)^F
	\end{equation*}
	and
	\begin{equation*}
	\operatorname{per}(aP_\theta + bP_\tau) = (a^l+b^l)^r(a+b)^F.
	\end{equation*}
\end{cor}

\begin{cor}
	Let $x_1, \dots , x_n$ be linearly independent vectors in $\mathbb{C}^n$ and $X$ be a matrix which $i$-th row is $x_i^T$. Denote $X^{-1} = (z_{ij})$. For $i=1, \dots , n$, denote $y_i = \sum\limits_{i=1}^n(az_{i \theta(j)}+bz_{i \tau(j)})e_i$, we have
	\begin{equation*}
	\langle x_1 \ast \cdots \ast x_n, y_1 \ast \cdots \ast y_n \rangle = \lvert G \rvert d_\chi^G(aP_\theta + bP_\tau).
	\end{equation*}
	\begin{proof}
		Let $A = (\langle x_i, y_j \rangle)$. Then $A = aP_\theta + bP_\tau$. By (\ref{eq21}), the corollary holds.
	\end{proof}
\end{cor}

In matrix theory, lots of properties are related to PSD matrices. It is worth to remark the following proposition.
\begin{prop}\label{PSD}
	Let $a,b \in \mathbb{C}$ and $\theta,\tau\in S_n$. Then $aP_\theta + bP_\tau$ is PSD if and only if $aP_\theta + bP_\tau = kI+mP_\pi$, where $k, m \in \mathbb{R}$ with $k \geq \lvert m \rvert$ and $\pi = \pi^{-1}$.
\end{prop}
\begin{proof}
	Suppose that $A = kI+mP_\pi$, where $k, m \in \mathbb{R}$ with $k \geq \lvert m \rvert$ and $\pi = \pi^{-1}$. For each $x = (x_1, \dots, x_n)^T \in \mathbb{C}^n$, we have
	\begin{equation*}
	x^{\ast}Ax = \sum\limits_{i=1}^n k\overline{x_i}x_i+m\overline{x_{\pi(i)}}x_i=\frac{1}{2}\sum\limits_{i=1}^n k(\overline{x_i}x_i+\overline{x_{\pi(i)}}x_{\pi(i)})+m(\overline{x_{\pi(i)}}x_i + \overline{x_i}x_{\pi(i)}).
	\end{equation*}
	Since $k \geq \lvert m \rvert$ and $\overline{x_i}x_i+\overline{x_{\pi(i)}x_{\pi(i)}} \geq \lvert\overline{x_{\pi(i)}}x_i + \overline{x_i}x_{\pi(i)}\rvert$ for each $i\in[n]$, we have $x^{\ast}Ax \geq 0$. 
	
	We will prove the converse by showing that $aP_\theta + bP_\tau$ is PSD only if it satisfies at least one of the following conditions:
	\begin{enumerate}
		\item $aP_\theta + bP_\tau = cI$ for some $c \geq 0$,
		\item $\theta = id$, $\tau^{-1}=\tau$ and $a \geq \lvert b \rvert$,
		\item $\theta^{-1} = \theta$, $\tau = id$ and $b \geq \lvert a \rvert$,
		\item $\theta, \tau$ are distinct nontrivial self-inverse with $[n] = \operatorname{Fix}(\theta) \cup \operatorname{Fix}(\tau)$ and $a = b \geq 0$.
	\end{enumerate}
Suppose that $aP_\theta + bP_\tau$ is PSD not satisfying (1)-(3). Then $\theta$ and $\tau$ are self-inverse. If $[n] \neq \operatorname{Fix}(\theta)\cup \operatorname{Fix}(\tau)$, there exists $i \not\in \operatorname{Fix}(\theta)\cup \operatorname{Fix}(\tau)$. This implies that
\begin{equation*}
(aP_\theta + bP_\tau)_{ii} = 0.
\end{equation*}
Note that if a diagonal entry of a PSD matrix is 0, then every entry in the row and column containing it is also 0, \cite[page 9-7]{LHogbenB}, which implies that $a = 0 = b$. Thus $[n] = \operatorname{Fix}(\theta) \cup \operatorname{Fix}(\tau)$. Since $\theta, \tau$ are not identities, $\operatorname{Fix}(\theta)\setminus \operatorname{Fix}(\tau)$ and $\operatorname{Fix}(\tau)\setminus \operatorname{Fix}(\theta)$ are nonempty. Let $i \in \operatorname{Fix}(\theta) \setminus \operatorname{Fix}(\tau)$ and $j = \tau(i)$. Then $(i,i)$-entry and $(i,j)$-entry are only two nonzero entries of $i$-th row. Similarly, $(j,i)$-entry and $(j,j)$-entry are only two nonzero entries of $j$-th row. Put
\begin{equation*}
\hbox{$x = e_i + e_j$ and $y = e_i - e_j$.}
\end{equation*}
Then
\begin{equation*}
\hbox{$2(a+b) = x^\ast Ax \geq 0$ and $2(a-b) = y^\ast A y \geq 0$.}
\end{equation*}
Thus $a \geq \lvert b \rvert$. Similarly, by replacing $y$ with $z = -y = e_j-e_i$, we have, $b \geq \lvert a \rvert$, and thus $a = b \geq 0$. Now, we have that $aP_\theta + bP_\tau$ satisfies at least one of the conditions. Since each condition implies that $aP_\theta + bP_\tau = kI+mP_\pi$, where $k,m\in \mathbb{R}$ with $k\geq \lvert m \rvert$ and $\pi = \pi^{-1}$. \end{proof}

We know that, for a character $\chi$ of a subgroup $G$ of $S_n$, $\chi(id) \geq \lvert\chi(\sigma)\rvert$ for each $\sigma \in S_n$. Let $N = kI +mP_\pi$, where $k,m\in \mathbb{R}$ with $k \geq \lvert m \rvert$ and $\pi = \pi^{-1}$. Then $N$ is PSD by Proposition \ref{PSD}. Let $\sigma \in X(id, \pi)$. If $\pi = C_1 \cdots C_r$, there exists $I \subseteq [r]$ such that $\sigma = \prod\limits_{i \in I}C_i$. We know that any self-inverse permutation is a product of disjoint transpositions. This means that $C_i$ is a transposition for each $i = 1, \dots, r$. Then $\sigma = \sigma^{-1}$. This implies that $\chi(\sigma) \in \mathbb{R}$ and $\prod\limits_{i=1}^n N_{i \, \sigma(i)}$ is a positive real number. Thus
\begin{eqnarray*}
\frac{1}{\chi(id)}d_\chi^G(N) & = & \sum\limits_{\sigma \in X(id, \pi)\cap G}\frac{\chi(\sigma)}{\chi(id)}\prod\limits_{i=1}^n N_{i \, \sigma(i)}\\
& \leq & \sum\limits_{\sigma \in X(id, \pi)\cap G}\frac{\lvert\chi(\sigma)\rvert}{\chi(id)}\prod\limits_{i=1}^n N_{i \, \sigma(i)}\\
& \leq &\sum\limits_{\sigma \in X(id, \pi)}\prod\limits_{i=1}^n N_{i \, \sigma(i)}\\
& = & \operatorname{per}(N).
\end{eqnarray*}
We now conclude that:
\begin{thm}The permanent dominance conjecture holds for a PSD matrix in the form $aP_\theta+bP_\tau$.
\end{thm}

In matrix theory, many inequalities are related to singular values. Let $A = aP_\theta + bP_\tau$, where $\theta, \tau \in S_n$. Here, we remark the form of a singular value of $A$. It is obvious that $A^{\ast} = \overline{a}P_{\theta^{-1}} + \overline{b}P_{\tau^{-1}}$. We have that
\begin{eqnarray*}
A^{\ast}A & = & (\lvert a \rvert^2 + \lvert b \rvert^2)I + \overline{a}bP_{\theta^{-1}\tau} + a\overline{b}P_{\tau^{-1}\theta}\\
& = & (\lvert a \rvert^2 + \lvert b \rvert^2)I + \overline{a}bP_{\theta^{-1}\tau} + (\overline{a}bP_{\theta^{-1}\tau})^{\ast}.
\end{eqnarray*}
Denote $X = \overline{a}bP_{\theta^{-1}\tau}$. Note that $X^{\ast}X = \lvert ab \rvert^2I = XX^{\ast}$. Then $X$ and $X^{\ast}$ are simultaneously diagonalizable. If $\lambda_1, \dots,\lambda_n$ are eigenvalues of $X$, then an eigenvalue of $X+X^{\ast}$ must be in the form
\begin{equation*}
\lambda_j + \overline{\lambda_{\sigma(j)}},
\end{equation*}
for some $\sigma \in S_n$. Furthermore, if $\theta^{-1}\tau$ has a cycle structure $[l_1, \dots, l_m]$, then $\lambda_j = \overline{a}b\xi$, where $\xi$ is $l_k$-th root of unity for some $k = 1, \dots, m$. Note that $X + X^{\ast}$ is hermitian. So, $\lambda_j + \overline{\lambda_{\sigma(j)}} \in \mathbb{R}$. We have that
\begin{equation*}
\lambda_j + \overline{\lambda_{\sigma(j)}} = \overline{a}be^{i\theta_1} + a\overline{b}e^{i\theta_2} \in \mathbb{R},
\end{equation*}
where $0 \leq \theta_1,\theta_2 < 2\pi$. This happens if and only if $\sin(\theta_1) = -\sin(\theta_2)$. Thus, an eigenvalue of $X + X^{\ast}$ is in the form $2Re(\lambda_j)$ which implies that $(\lvert a \rvert^2 + \lvert b \rvert^2 + 2Re(\lambda_j))^{\frac{1}{2}}$ is a singular value of $A^{\ast}A$. The following proposition is received by applying (\ref{sveq}).
\begin{prop}
Let $\chi$ be a linear character of a subgroup $G$ of $S_n$. For any $\theta,\tau \in S_n$ and $a,b\in \mathbb{C}$,
\begin{equation*}
\lvert \sum\limits_{\sigma\in X(\theta,\tau) \cap G}\overline{\chi}(\sigma)a^{n - t_\sigma - F} b^{t_\sigma}(a+b)^F \rvert^2 \leq \frac{1}{n}\sum\limits_{j=1}^n(\lvert a \rvert^2 + \lvert b \rvert^2 + 2Re(\lambda_j))^{n},
\end{equation*}
where $\lambda_1, \dots, \lambda_n$ are the eigenvalues of $\overline{a}bP_{\theta^{-1}\tau}$.	
\end{prop}
This proposition is an example of matrix inequality relating to singular values and generalized matrix functions. Many inequalities can be found in \cite{RMerrisB,XZhanB,FZhangB}.
\section{generalized matrix functions on cousins of permutation matrices}

In this section, by applying some technique in the previous section, we generalized the result to the case of block matrices.

For a permutation $f$ on $[n]$ and a nonnegative integers $x, y$, we define a function $f_{x,y}:x+[n] \rightarrow y+[n]$ by
\begin{equation*}
\hbox{$f_{x,y}(x+i) = y+f(i)$ for each $i \in [n]$.}
\end{equation*}

For a function $f:X\rightarrow Y$ and $g:W\rightarrow Z$, which $X$ and $W$ are disjoint, we denote $f\sqcup g: X \cup W \rightarrow Y \cup Z$ to be the extension of $f$ and $g$.

Let $a_1, \dots, a_n, b_1, \dots ,b_n \in \mathbb{C}$ and $\theta, \tau \in S_n$. Define $n \times n$ block matrices $P, Q$ whose each block is an $m \times m$ complex matrix by
\begin{equation*}
[P]_{i \, j} = \left\{
\begin{array}{ll}
a_i P_{\theta_i}, & \hbox{if $j = \theta(i)$,} \\
0, & \hbox{otherwise,}
\end{array}
\right.
\end{equation*}
and
\begin{equation*}
[Q]_{i \, j} = \left\{
\begin{array}{ll}
b_i P_{\tau_i}, & \hbox{if $j = \tau(i)$,} \\
0, & \hbox{otherwise,}
\end{array}
\right.
\end{equation*}
where $\theta_i$ and $\tau_i$ are permutations in $S_m$ for each $i=1,\dots,n$. Here, $[P]_{i \, j}$ and $P_{i \, j}$ denote $(i, j)$-block and $(i, j)$-entry of $P$, respectively. Let $M = P+Q$. This matrix is one of the ways to consider a generalized permutation matrix as a block matrix. If $n = 1$, $M$ is a linear sum of two permutation matrices. We will evaluate the value of $d_\chi^G(M)$. Let $\alpha = \bigsqcup\limits_{i=1}^m \theta_{i-1, \theta(i) - 1}$ and $\beta = \bigsqcup\limits_{i=1}^m \tau_{i-1, \tau(i) - 1}$. Then $\alpha, \beta \in S_{mn}$. Now, suppose that $\alpha^{-1}\beta = C_1 \cdots C_r$ is the disjoint cycles expression.  Because $M_{i \, \sigma(i)} = 0$ if there exists $i\in [n]$ such that $\sigma(i)$ is not $\theta^{-1}(i)$ nor $\tau^{-1}(i)$, we obtain that the index of the summand of $d_\chi^G(M)$ is $X(\alpha,\beta) \cap G$ for each subgroup $G$ of $S_{mn}$.
For each $\sigma \in X(\alpha,\beta)$, denote a subset $I$ of $[r]$ with $\sigma = \alpha\prod\limits_{i\in I} C_i$ by $I^{\sigma}$. Suppose that $l_i$ is the length of $C_i$. We denote $l_{i, j}$ the number of element of $\{jm, jm-1, \dots, jm-m\}$ containing in $C_i$. We also denote $F_j$ the number of element of $\{jm, jm-1, \dots, jm-m\}$ which is a fixed point of $\alpha^{-1}\beta$.
The following theorem is true.
\begin{thm}\label{thm4}
	Let $\chi$ be a character of subgroup $G$ of $S_{mn}$. Define an $mn \times mn$ matrix $M$ as the above discussion. Then
	\begin{equation*}
	d_\chi^G(M) = \prod\limits_{i=1}^n(a_i + b_i)^{F_i}\sum\limits_{\sigma\in X(\alpha, \beta)\cap G}\overline{\chi}(\sigma)\prod\limits_{i \not\in I^{\sigma}}\prod\limits_{j=1}^na_j^{l_{i, j}} \prod\limits_{i \in I^{\sigma}}\prod\limits_{j=1}^nb_j^{l_{i, j}}.
	\end{equation*}
\end{thm}
\begin{proof}
	By the similar reasoning to the discussion preceding Theorem \ref{form1}, we have
	\begin{eqnarray*}
	\prod\limits_{i=1}^{mn}M_{i \, \sigma(i)} & = & \sum\limits_{I \subseteq [mn]}\prod\limits_{i \in I}P_{i\, \sigma(i)}\prod\limits_{i \not\in I}Q_{i\, \sigma(i)}\\
	& = & \sum\limits_{A\subseteq \operatorname{Fix}(\alpha^{-1}\beta)}\prod\limits_{i \in I^\sigma_A}P_{i\, \sigma(i)}\prod\limits_{i \not\in I^{\sigma}_A}Q_{i\, \sigma(i)}\\
	& = & \sum\limits_{A\subseteq \operatorname{Fix}(\alpha^{-1}\beta)}\prod\limits_{j=1}^n\prod\limits_{i \in (I^\sigma_A)_j}a_j(P_{\theta_j})_{i\, \sigma(i)}\prod\limits_{i \not\in (I^{\sigma}_A)_j}b_j(P_{\tau_j})_{i\, \sigma(i)},
	\end{eqnarray*}
where $(I^{\sigma}_A)_j = I^{\sigma}_A \cap ((j-1)m,jm]$.  By the same way as we did with $t_\sigma$ and $F$ in (\ref{tF}) in the discussion preceding Theorem \ref{form1}, we have
\begin{equation*}
\prod\limits_{i=1}^{mn}M_{i \, \sigma(i)} = \prod\limits_{i \not\in I^{\sigma}}\prod\limits_{j=1}^na_j^{l_{i, j}} \prod\limits_{i \in I^{\sigma}}\prod\limits_{j=1}^nb_j^{l_{i, j}},
\end{equation*}
which completes the proof.
\end{proof}

Like the index sum of the formula for $d_\chi^G(aP_\theta + bP_\tau)$, the index of summand in the formula of $d_\chi^G(M)$ is reduced from $G$ to $X(\alpha, \beta) \cap G$. The value of $\prod\limits_{i \not\in I^{\sigma}}\prod\limits_{j=1}^na_j^{l_{i, j}} \prod\limits_{i \in I^{\sigma}}\prod\limits_{j=1}^nb_j^{l_{i, j}}$ can be calculated directly by counting the number of elements of $\{jm, jm-1, \dots jm-m\}$ containing in each cycle of $\alpha^{-1}\beta$.

\begin{exaple}
	Let
	\begin{equation*}
	M = \begin{bmatrix}
	0 & 0 & -i & 0    & 0 & -2 & 0 & 0\\
	0 & -i & 0 & 0    & 0 & 0 & -2 & 0\\
	0 & 0 & 0 & -i    & -2 & 0 & 0 & 0\\
	-i & 0 & 0 & 0    & 0 & 0 & 0 & -2\\
	
	3 & 0 & 0 & 0    & 0 & 0 & 0 & 2\\
	0 & 3 & 0 & 0    & 0 & 0 & 2 & 0\\
	0 & 0 & 3 & 0    & 0 & 2 & 0 & 0\\
	0 & 0 & 0 & 3    & 2 & 0 & 0 & 0
	\end{bmatrix}.
	\end{equation*}
	Here, we can consider $M$ as $\begin{bmatrix}
	-iP_{(143)} & -2P_{(132)}\\
	3I & 2P_{(14)(23)}
	\end{bmatrix} = \begin{bmatrix}
	-iP_{(143)} & 0\\
	0 & 2P_{(14)(23)}
	\end{bmatrix} + \begin{bmatrix}
	0 & -2P_{(132)}\\
	3I & 0
	\end{bmatrix}$ (heart form). Comparing to Theorem \ref{thm4}, we have $m = 4, n = 2,\theta_{1} = (143), \theta_2 = (132), \theta_3 = id, \theta_4 = (14)(23)$. We derive that
	\begin{equation*}
	\hbox{$\alpha = \theta_1 \sqcup \theta_4 = (143)(58)(67)$ and $\beta = \theta_2 \sqcup \theta_3 = (153726)(48)$.}
	\end{equation*}
	So, $\alpha^{-1}\beta = (18)(27)(36)(45)$ which implies that $F_i = 0, l_i = 2$ and $l_{i,j} = 1$ for all $i=1,2,3,4$ and all $j =1,2$. In this example, $a_1 = -i, a_2 = 2, b_1 = -2, b_2 = 3$. By Theorem \ref{thm4},
	\begin{equation*}
	\hbox{$\operatorname{per}(M) = (-6)^4+4(-2i)(-6)^3+6(-2i)^2(-6)^2+4(-2i)^3(-6)+(-2i)^4= 448+1536i$}
	\end{equation*}
	and
	\begin{equation*}
	\hbox{$\det(M) = (-6)^4-4(-2i)(-6)^3+6(-2i)^2(-6)^2-4(-2i)^3(-6)+(-2i)^4= 448-1536i$.}	
	\end{equation*}
\end{exaple}

In \cite{RS2019}, we introduced an $n \times n$ symmetric matrix $S_\theta$, where $\theta \in S_n$. The formula for a generalized matrix function of $S_\theta$ were also investigated. This matrix was used to prove a necessary and sufficient condition for the equality of two generalized matrix function on the set of all symmetric matrices. Here, we recall the definition of $S_\theta$. For $\theta \in S_n$, $S_\theta$ is a symmetric matrix with
\begin{equation*}
(S_\theta)_{i \, j} = \left\{
\begin{array}{ll}
1, & \hbox{if $\theta(i)=j$ or $\theta^{-1}(i)=j$,} \\
0, & \hbox{otherwise.}
\end{array}
\right.
\end{equation*}

A matrix $S_\theta$ can be considered as a special case of the matrix $M$. This fact is obtained by putting $n = 1$, $\alpha = \theta$, $\beta = \theta^{-1}$, $a_i = 1 = b_i$ if $\theta^2(i) \neq i$, otherwise $a_i = \frac{\displaystyle 1}{\displaystyle 2} = b_i$. Moreover, by the definition of $S_\theta$, we immediately have that
\begin{equation}\label{eq000}
d_\chi^G(P_\theta + P_{\theta^{-1}}) = 2^{F+2t}d_\chi^G(S_\theta),
\end{equation}
where $F$ and $t$ are numbers of fixed point and transpositions of $\theta$, respectively.

For $\theta,\tau \in S_n$, we have that $P_\theta P_\tau = P_{\theta\tau}$. This property does not hold for $S_\theta$ and $S_\tau$ in general, but it is also true for some pairs of permutations as we can see in the following proposition.
\begin{prop}\label{thm02}
	Let $\theta,\tau \in S_n$ with each cycle from the disjoint expression of $\theta$ and each from the disjoint expression of $\tau$ are disjoint.
	Then $S_\theta S_\tau = S_{\theta\tau}$.
\end{prop}
\begin{proof}
	We know that
	\begin{equation*}
	(S_\theta S_\tau)_{i \, j} = \sum\limits_{k = 1}^n (S_\theta)_{i \, k}(S_\tau)_{k \, j}.
	\end{equation*}
	If $i \in \operatorname{Fix}(\theta)$, then $\theta\tau(i) = \tau\theta(i) = \tau(i)$.
	So, $(S_\theta S_\tau)_{i \, j} = (S_\tau)_{i \, j} = (S_{\theta\tau})_{i \, j}$.
	Suppose that $i \in \operatorname{Fix}(\theta)^c$.
	Since $\theta$ and $\tau$ are disjoint, $\theta(i), \theta^{-1}(i) \not\in \operatorname{Fix}(\tau)^c$.
	If $\theta(i) = \theta^{-1}(i)$, then $(S_\theta S_\tau)_{i \, j} = (S_\theta]_{i \, \theta(i)}(S_\tau)_{\theta(i) \, j}$.
	This implies that
	\begin{equation*}
	(S_{\theta\tau})_{i \, j} = \left\{
	\begin{array}{ll}
	1, & \hbox{if $j = \theta(i)$,} \\
	0, & \hbox{otherwise.}
	\end{array}
	\right.
	\end{equation*}
	Assume that $\theta(i) \neq \theta^{-1}(i)$. Then
	\begin{equation*}
	(S_\theta S_\tau)_{i \, j} = (S_\theta)_{i \, \theta(i)}(S_\tau)_{\theta(i) \, j} + (S_\theta)_{i \, \theta^{-1}(i)}(S_\tau)_{\theta^{-1}(i) \, j}.
	\end{equation*}
	Thus
	\begin{equation*}
	(S_{\theta\tau})_{i \, j} = \left\{
	\begin{array}{ll}
	1, & \hbox{if $j = \theta(i)$ or $j = \theta^{-1}(i)$,} \\
	0, & \hbox{otherwise.}
	\end{array}
	\right.
	\end{equation*}
	Because $\theta(i) \in \operatorname{Fix}(\tau)$, we have $\theta(i) = \tau\theta(i) = \theta\tau(i)$. We obtain that $(S_\theta S_\tau)_{i \, j} = (S_\theta)_{i \, j} = (S_{\theta\tau})_{i \, j}$. Hence $S_\theta S_\tau = S_{\theta\tau}$.
\end{proof}

Let $\pi$ be a cycle in $S_n$. Then $\pi^2$ is a cycle if and only if an order of $\pi$ is odd. If $\pi$ is a $2k$-cycle, then $\pi^2$ is a product of two $k$-cycles. Moreover, $k$ is odd if and only if these two cycles are even. By Theorem \ref{form1} and (\ref{eq000}), this implies that
\begin{equation*}
\det(S_\pi) = \left\{
\begin{array}{ll}
-1, & \hbox{if an order of $\pi$ is 2,} \\
2, & \hbox{if an order of $\pi$ is odd,} \\
-4, & \hbox{if an order of $\pi$ is $4k+2$ for some $k \in \mathbb{N}$,}\\
0, & \hbox{if an order of $\pi$ is $4k$ for some $k \in \mathbb{N}$.}
\end{array}
\right.
\end{equation*}
By applying Proposition \ref{thm02}, the following result holds.
\begin{cor}
	Let $\theta \in S_n$. Then the decomposition of $\theta$ contains a $4m$-cycle if and only if $\det(S_\theta) = 0$. Moreover,
	\begin{equation*}
	\det(P_\theta + P_{\theta^{-1}}) = (-1)^{s+t}2^{F+r+2(s+t)},
	\end{equation*}
	and
	\begin{equation*}
	\det(S_\theta) = (-1)^{s+t}2^{r+2s},
	\end{equation*}
	where $F,r,s,t$ are numbers of fixed point, cycles of odd order, order $4k+2$ for some $k \in \mathbb{N}$, and order $2$ of $\theta$, respectively.
\end{cor}
\section*{Acknowledgements}
The second author would like to thank Department of Mathematics and Faculty of Science, Naresuan University for the financial support.

\end{document}